\documentclass[11pt]{article}

\usepackage{graphicx}
\usepackage{amsmath}
\usepackage{amsfonts}
\usepackage{fancyhdr} 

\newtheorem{lem}{Lemma}
\newtheorem{thm}{Theorem}

\newtheorem{cor}[lem]{Corollary}

\newcommand{\fff}{\mathcal{F}}
\newcommand{\sss}{\mathcal{S}}

\newcommand{\real}{\mathbb{R}}

\newcommand{\poly}{\mathcal{P}}
\newcommand{\ooo}{\mathcal{O}}
\newcommand{\pr}{^{\prime}}
\newcommand{\ddd}{\mathbb{D}}

\newcommand{\complex}{\mathbb{C}}

\newcommand{\hright}{H^2_{\complex_+}}
\newcommand{\hdisk}{H^2_{\ddd}}

\newcommand{\mmm}{\mathcal{M}}
\newcommand{\tmm}{\widetilde{\mathcal{M}}}
\newcommand{\lll}{\mathcal{L}}
\newcommand{\hhh}{\mathcal{H}}
\newcommand{\ltwo}{L^2(\real_+)}
\newcommand{\lltwo}{l^2(\mathbb{Z}_+)}

\newcommand{\aaa}{\mathcal{A}}

\newenvironment{pf}{\noindent {\em Proof}.\ \ }{\hspace*{\fill}\rule{.5ex}{1.4ex}\,}

\title{Identification of minimum phase preserving operators on the half line}
\author{Peter C.~Gibson\footnote{Corresponding author} \footnote{Dept.~of Mathematics \& Statistics, York University, 4700 Keele St., Toronto, Ontario, Canada, M3J~1P3, $\mathtt{pcgibson@yorku.ca}$} and Michael P.~Lamoureux\footnote{Dept.~of Mathematics \& Statistics, University of Calgary, 2500 University Dr. NW, Calgary, Alberta, Canada, T2N~1N4, $\mathtt{mikel@ucalgary.ca}$}}

\begin{document}

\maketitle
\begin{abstract}
Minimum phase functions are fundamental in a range of applications, including control theory, communication theory and signal processing.   A basic mathematical challenge that arises in the context of geophysical imaging is to understand the structure of linear operators preserving the class of minimum phase functions.  The heart of the matter is an inverse problem: to reconstruct an unknown minimum phase preserving operator from its value on a limited set of test functions.   This entails, as a preliminary step, ascertaining sets of test functions that determine the operator, as well as the derivation of a corresponding reconstruction scheme.  In the present paper we exploit a recent breakthrough in the theory of stable polynomials to solve the stated inverse problem completely.   We prove that a minimum phase preserving operator on the half line can be reconstructed from data consisting of its value on precisely two test functions.  And we derive an explicit integral representation of the unknown operator in terms of this data.  A remarkable corollary of the solution is that if a linear minimum phase preserving operator has rank at least two, then it is necessarily injective. 
\end{abstract}

\newpage

%\tableofcontents

\section{Introduction\label{sec-Introduction}}

\pagestyle{fancyplain}

In the space $\ltwo$, which is the basic model for causal signals of a continuous variable, the term minimum phase refers to functions $f$ which are determined by their power spectrum, i.e. the modulus $|\mathcal{F}f|$ of the Fourier transform, and which among all functions having a given power spectrum maximize the partial energies $\int_0^T|f|^2$, for all $T>0$.  Minimum phase functions are fundamental in myriad contexts, including control theory, communication theory and signal processing.  This paper presents a full solution to the following inverse problem, which is motivated by geophysical considerations.  Let $A$ be an unknown linear operator on $\ltwo$ with the property that the function $Af$ is the translate of a minimum phase function whenever $f$ is the translate of a minimum phase function.  Which sets of functions $\mathcal{S}\subset\ltwo$ determine $A$?  And how can $A$ be reconstructed from its values $A\mathcal{S}$ on such sets?  The main result of the paper, Theorem~\ref{thm-Reformulation} in Section~\ref{sec-Identification}, gives an explicit formulation for such an operator $A$ in terms of its value at the two functions, $e^{-t}(1-t)$ and $e^{-t}t$.  The proof uses a new result in the theory of stable polynomials that characterizes all linear operators mapping polynomials having zeros outside the open unit disk to zero-free analytic functions on the disk.  Hardy spaces on the half plane and disk serve as the basic framework for our arguments.    

In geophysics there is a long standing idea---see \cite{ShTr:1965}, for example---that in a horizontally stratified absorptive earth with vertically traveling plane compressional waves (represented as functions of travel time), transmitted waves are translates of minimum phase functions. And there have been field experiments, \cite{ZiBo:1993}, supporting the assertion that, for instance, the source signature of a dynamite blast is minimum phase.  This provides the conceptual basis for the inverse problem considered in the present article, supposing that the aforementioned minimum phase waveforms are the output of a linear operator that appropriately encodes the material properties of the earth. In the present article no attempt is made to analyze the physical model of a layered earth, nor to explicitly represent solution operators for the attendant governing equations.  Rather, the minimum phase hypothesis is taken as a starting point, and the aim is to investigate its consequences.  However, the precise results obtained in the present paper offer the prospect of practical significance for geophysical imaging, as well as implications for the wider arena in which minimum phase signals play an important role. 

The paper is organized as follows.   Section~\ref{sec-Notation}, below, assembles the necessary notation and background material required to solve the main inverse problem.  The background material is divided into four subsections:  Section~\ref{sec-Hardy} covers the notation and standard factorization theorems for the Hardy-Hilbert spaces on the half plane and disk; Section~\ref{sec-Isometry} establishes an isometric isomorphism between $\ltwo$ and the Hardy space on the disk, which helps to streamline later definitions and proofs; Section~\ref{sec-Minimum} covers the technicalities of minimum phase, as it relates to both the discrete and continuous settings; and Section~\ref{sec-Stable} quotes two results from the theory of stable polynomials that are needed to prove the main results.   

There are two main steps toward solving the main inverse problem.  The first step is to derive a structural characterization of linear operators on $\ltwo$ that preserve translated minimum phase operators. This is carried out in Section~\ref{sec-Characterization}, where Theorem~\ref{thm-Characterization} describes minimum phase preserving operators in terms of product-composition operators on Hardy space on the disk, and Theorem~\ref{thm-Characterization2} describes them in terms of product-composition operators on Hardy space on the half plane.  

The final step toward solving the main inverse problem is to reconstruct an unknown minimum phase operator in terms of test data.  This is done in Section~\ref{sec-Identification}, where two explicit ways of identifying an unknown minimum phase preserving operator are given, in Theorems~\ref{thm-Identification} and \ref{thm-Reformulation} respectively.  The characterization from the previous section is incorporated into the statement of Theorem~\ref{thm-Reformulation}, making the latter the main result of the paper.  

In Section~\ref{sec-Remarks}, the case of untranslated minimum phase signals is handled briefly in Section~\ref{sec-Operators}.  Then Section~\ref{sec-Hypothesis} points out a remarkable corollary of the main result, that a linear operator that preserves minimum phase is necessarily injective, with the exception of the rank 1 case.  

Lastly, Section~\ref{sec-Conclusion} briefly concludes the paper.

\section{Notation and background results\label{sec-Notation}}

\subsection{The Hardy-Hilbert spaces on the half plane and disk\label{sec-Hardy}}

Let $\hright$ denote the Hardy space consisting of all analytic functions on the open right half plane $\complex_+$ 
whose square integrals along lines of constant positive real part are uniformly bounded, 
and let $\hdisk$ denote the Hardy space of all analytic functions on the open unit disk $\ddd$ whose Taylor coefficients are square summable.  Both $\hright$ and $\hdisk$ are separable Hilbert spaces.   Indeed the operator
\[
\Phi:\hright\rightarrow\hdisk
\]
defined by 
\[
(\Phi F)(z)=\frac{2\sqrt{\pi}}{1+z}\,F\!\left(\frac{1-z}{1+z}\right)
\]
for every $F\in\hright$ and $z\in\ddd$ is a surjective isometry, with inverse
\[
(\Phi^{-1}G)(w)=\frac{1}{\sqrt{\pi}(1+w)}\,G\!\left(\frac{1-w}{1+w}\right).
\]
See \cite[Chapter~8]{Ho:1962}.  

A function $G\in\hdisk$ is defined to be outer if $G$ is not the zero function and there exists a scalar $\lambda$ of modulus 1 such that for every $z\in\ddd$
\[
G(z)=\lambda\exp\!\left(\frac{1}{2\pi}\int_0^{2\pi}\frac{e^{i\theta}+z}{e^{i\theta}-z}\log\left|G(e^{i\theta})\right|\,d\theta\right);
\]
$G$ is defined to be inner if $|G(e^{i\theta})|=1$ for almost every point $e^{i\theta}$ belonging to the unit circle $S^1$.   It is a fundamental fact that, apart from the zero function, every $G\in\hdisk$ has a factorization of the form 
\[
G=G_{\rm inner}G_{\rm outer},
\]
where $G_{\rm inner}\in\hdisk$ is inner and $G_{\rm outer}\in\hdisk$ is outer; this factorization is unique up to multiplication by a scalar in $S^1$.   

A non-zero function $F\in\hright$ is defined to be outer if there exists a scalar $\lambda\in S^1$ such that for every $w\in\complex_+$,
\[
F(w)=\lambda\exp\!\left(\frac{1}{\pi}\int_{-\infty}^\infty\frac{yw+i}{y+iw}\log\left|F(iy)\right|\frac{dy}{1+y^2}\right).
\]
Using the fact that $1+z\in\hdisk$ is outer, direct calculation shows that a function $F\in\hright$ is outer if and only if its image $\Phi F\in\hdisk$ is outer.  

The inner-outer factorization in $\hdisk$ translates via $\Phi$ to $\hright$ as follows.  Given $F\in\hright$ let $G=G_{\rm inner}G_{\rm outer}=\Phi F$, and let $m:\complex\cup\{\infty\}\rightarrow\complex\cup\{\infty\}$ denote the M\"obius transformation
\begin{equation}\label{Mobius}
m(z)=m^{-1}(z)=\frac{1-z}{1+z}, 
\end{equation}
which maps $\ddd$ onto $\complex_+$ and vice versa. (In signal processing $m(z^{-1})$ is known as the bilinear transform---see \cite[Chapter~7]{OpSc:2010}.)  Observe that 
\[
\Phi^{-1}(G_{\rm inner}G_{\rm outer})=(G_{\rm inner}\circ m)\!\cdot\!(\Phi^{-1}G_{\rm outer}).
\]
This induces the factorization $F=F_{\rm inner}F_{\rm outer}$,  where the factor
\[
F_{\rm inner}=G_{\rm inner}\circ m
\]
has modulus 1 on the imaginary line $w=iy$, and the factor
\[
F_{\rm outer}=\Phi^{-1}G_{\rm outer}
\]
is outer.  Note that the factorization on $\hright$ differs from the one on $\hdisk$ in that the factor $F_{\rm inner}$, which is not square integrable on the imaginary line, does not belong to $\hright$.

\subsection{The isometry $\hhh$ mapping $\ltwo$ onto $\hdisk$\label{sec-Isometry}}

The Fourier-Laplace transform 
\[
\lll :L^2(\real_+)\rightarrow\hright
\]
defined by 
\[
(\lll f)(w)=\frac{1}{\sqrt{2\pi}}\int_0^\infty f(t)e^{-wt}\,dt
\]
for every $f\in L^2(\real_+)$ and $w\in\complex_+$ is a surjective isometry.  Its boundary function is the Fourier transform of $f$, denoted $\mathcal{F}f$; that is,
\[
(\lll f)(iy)=(\mathcal{F}f)(y)\rule{10pt}{0pt}\mbox{ for a.e. }y\in\real.
\]  
Let $\hhh:L^2(\real_+)\rightarrow\hdisk$ denote the composition 
\[
\hhh=\Phi\lll.
\]
Note that $\hhh$ is a surjective isometry, since each of its factors is.   The map $\hhh$ is given explicitly by the formula
\[
(\hhh f)(z)=\frac{\sqrt{2}}{1+z}\int_0^\infty f(t)\exp\!\left(\frac{t(z-1)}{z+1}\right)\,dt
\]
for every $f\in\ltwo$ and $z\in\ddd$.   Its inverse mapping   
\[
\hhh^{-1}:\hdisk\rightarrow\ltwo
\]
is given by the formula
\[
(\hhh^{-1}G)(t)=\frac{1}{\pi\sqrt{2}}\int_{-\infty}^\infty G\!\left(\frac{1-iy}{1+iy}\right)\,\frac{e^{iyt}}{1+iy}\,dy
\]
for every $t\in\real_+$, provided the boundary function of $\Phi^{-1}G$ is in $L^1(\real)$; otherwise $\hhh^{-1}G$ may be computed as an oscillatory integral as per the standard theory for Fourier transforms on $L^2(\real)$.   

For each nonnegative integer $n$, let $\mu_n\in\hdisk$ denote the $n$th moment function
\[
\mu_n(z)=z^n.
\]
Since the functions $\mu_n$ comprise an orthonormal basis for $\hdisk$, and $\hhh^{-1}$ is an isometry, it follows that the functions $\hhh^{-1}\mu_n$ form an orthonormal basis for $\ltwo$.   These may be computed explicitly using the above formula for $\hhh^{-1}$, which, after some work, yields
\[
(\hhh^{-1}\mu_n)(t)=(-1)^n\sqrt{2}e^{-t}L_n(2t),
\]
for each $n\geq 0$, where $L_n$ is the $n$th Laguerre polynomial defined as 
\[
L_n(t)=\frac{e^t}{n!}\frac{d^n}{dt^n}\left(t^ne^{-t}\right).  
\]
In particular,  
\[
(\hhh^{-1}\mu_0)(t)=\sqrt{2}e^{-t} \mbox{ and } (\hhh^{-1}\mu_1)(t)=\sqrt{2}e^{-t}(2t-1).
\]

Let 
\[
\mathfrak{d}:\ltwo\rightarrow\lltwo
\]
 denote the map defined by the formula
\[
\mathfrak{d}f=\{a_n\}_{n=0}^\infty, 
\]
where for each $n\geq0$,
\[
 a_n=\int_0^\infty f(t)\overline{(\hhh^{-1}\mu_n)(t)}\,dt=
(-1)^n\sqrt{2}\int_0^\infty f(t)L_n(2t)e^{-t}\,dt.
\]
Note that, since the sequence of Fourier coefficients $\{c_n\}_{n=0}^\infty=\mathbf{F}G$ of a function $G\in\hdisk$ is given by the formula 
\[
c_n=\langle G,\mu_n\rangle=\frac{1}{2\pi}\int_0^{2\pi}G(e^{i\theta})e^{-in\theta}\,d\theta=\int_0^\infty (\hhh^{-1}G)(t)\overline{(\hhh^{-1}\mu_n)(t)}\,dt,
\]
the map $\mathfrak{d}$ simply associates the sequence of Fourier coefficients $\{a_n\}_{n=0}^\infty=\mathbf{F}\hhh f$ to a given continuous signal $f$.  Since the transform 
\[
\mathbf{F}:\hdisk\rightarrow\lltwo
\]
is an isometry, $\mathfrak{d}$ is therefore an isometry too.

\subsection{Minimum phase functions\label{sec-Minimum}}

The inverse of the transform $\mathbf{F}:\hdisk\rightarrow\lltwo$ is the $z$-transform 
\[
\mathbf{Z}:\lltwo\rightarrow\hdisk
\]
defined by the formula
\[
\left(\mathbf{Z}a\right)(z)=\sum_{n=0}^\infty a_n z^n 
\]
for every sequence $a=\{a_n\}_{n=0}^\infty\in\lltwo$ and $z\in\ddd$. A causal discrete signal in $\lltwo$ is defined to be minimum phase if its $z$-transform 
is an outer function, or in other words, if the given discrete signal is the sequence of Fourier coefficients of an outer function.  A minimum phase signal $\{a_n\}_{n=0}^\infty$ has maximally front loaded energy, in the sense that if $\{b_n\}_{n=0}^\infty\in\lltwo$ satisfies
\[
\left|\sum_{n=0}^\infty b_n z^n\right|=\left|\sum_{n=0}^\infty a_n z^n\right|
\]
for every $z\in\ddd$, then for every $N\geq 0$,
\[
\sum_{n=0}^N|b_n|^2\leq\sum_{n=0}^N|a_n|^2.
\]
See \cite{OpSc:2010}. 

A function $f\in\ltwo$ is defined to be minimum phase if its image $\hhh f\in\hdisk$ is outer.  As in the discrete case, if $f$ is minimum phase then its energy is maximally front loaded as follows.  
For every $g\in\ltwo$ such that $|\fff g|=|\fff f|$, and for every $T>0$,
\[
\int_0^T\left|f(t)\right|^2\,dt\geq\int_0^T\left|g(t)\right|^2\,dt.
\]
In the continuous setting the connection between minimum phase and the energy being maximally front loaded is much more delicate than in the discrete setting, where the two properties are well known to be equivalent.  Non-uniqueness of solutions to the Pauli problem---see\cite{Is:1996}---bedevils the continuous case; this will be treated in a separate paper.  
%See \cite{GiLa:Min2011} for a proof.   

The symbol $\mmm$ denotes the set of all minimum phase functions in $\ltwo$.  Given $f\in\mmm$ and $\tau\geq0$, the translate 
\[
(T_\tau f)(t)=\left\{\begin{array}{cc} f(t-\tau) &\mbox{ if }t\geq\tau\\
0&\mbox{ if }t<\tau
\end{array}\right.
\]
is minimum phase only if $\tau=0$.   Nevertheless, from the signal processing point of view, translates of minimum phase functions are of similar interest to minimum phase signals; they are also front loaded, but with respect to a given delay $\tau\geq0$.  The set of all translates of minimum phase signals is denoted $\tmm$.  That is, 
\[
\tmm=\left\{ T_\tau f\,\left|\,f\in\mmm\mbox{ and }\tau\geq0\right.\right\}.
\]

The foregoing definitions of minimum phase imply the following. 
\begin{thm}\label{thm-d}
As in Section~\ref{sec-Isometry}, let
\[
\mathfrak{d}:\ltwo\rightarrow\lltwo
\]
be defined by the formula
\[
\mathfrak{d}f=\{a_n\}_{n=0}^\infty, 
\]
where for each $n\geq0$,
\[
 a_n=\int_0^\infty f(t)\overline{(\hhh^{-1}\mu_n)(t)}\,dt=
(-1)^n\sqrt{2}\int_0^\infty f(t)L_n(2t)e^{-t}\,dt.
\]
Then $\mathfrak{d}$ is a surjective isometry with the property that $\mathfrak{d}f$ is minimum phase if and only if $f$ is minimum phase. 
\end{thm}
\begin{pf}
As noted at the end of Section~\ref{sec-Isometry}, $\mathfrak{d}=\mathbf{F}\hhh$ is an isometry since each of $\hhh$ and $\mathbf{F}$ is.  By definition, a function $f\in\ltwo$ is minimum phase if and only if $\hhh f$ is an outer function, which is true if and only if its sequence of Fourier coefficients 
\[
\mathfrak{d}f=\mathbf{F}\hhh f
\]
is minimum phase.  \end{pf}

\subsection{Stable polynomials and stable preserving operators\label{sec-Stable}}

The arguments given in the present paper use results from \cite{GiLa:Polya2011} concerning stable polynomials, necessitating some notation as follows.  Given a set $\Omega\subset\complex$, the algebra of all analytic functions $f:\Omega\rightarrow\complex$ is denoted by $\aaa(\Omega)$.   Given sets $\Omega\subset\Omega\pr\subset\complex$, a function $f:\Omega\pr\rightarrow\complex$ is said to be $\Omega$-stable if $f(z)\neq0$ for every $z\in\Omega$.  The set of all $\Omega$-stable functions in $\aaa(\Omega\pr)$ is denoted by $\sss(\Omega)$.   
The symbol $\poly(\Omega)$ denotes the set of all $\Omega$-stable polynomials.   Given $\psi\in\aaa(\Omega)$, 
\[
M_\psi:\aaa(\Omega)\rightarrow\aaa(\Omega)
\]
denotes the operation of multiplication by $\psi$,
\[
(M_\psi f)(z)=\psi(z)f(z) \rule{5pt}{0pt} \]
for every $f\in\aaa(\Omega)$ and $z\in\Omega$.  
Given an analytic function $\varphi:\Omega\rightarrow\Omega$, 
\[
C_\varphi:\aaa(\Omega)\rightarrow\aaa(\Omega)
\]
denotes the operation of (right) composition with $\varphi$,
\[
(C_\varphi f)(z)=f\circ\varphi(z)=f(\varphi(z))\]
for every $f\in\aaa(\Omega)$ and $z\in\Omega$.  The symbol $\complex[z]$ denotes the set of all univariate polynomials over $\complex$.  

\begin{thm}[from {\cite[Theorem~2]{GiLa:Polya2011}}]\label{thm-Structure}
Let $\Omega\subset\ddd$ be a non-empty connected open set.   If a linear map $A:\complex[z]\rightarrow\aaa(\ddd)$ has the property that 
\[
A(\poly(\ddd))\subset \sss(\Omega)\cup\{0\},
\]
then either:
\begin{enumerate}
\item there exist a function $\psi\in\sss(\Omega)$ and a linear functional $\nu:\complex[z]\rightarrow\complex$ such that $A(f)=\nu(f)\psi$, for all $f\in\complex[z]$; or 
\item there exist a function $\psi\in\sss(\Omega)$ and a non-constant function $\varphi\in\aaa(\ddd)$, where $\varphi(\Omega)\subset\ddd$, such that $A=M_\psi C_\varphi$.  
\end{enumerate}
\end{thm}

Note that there is no topology attributed to the vector space $\complex[z]$, and hence no hypothesis of continuity in the above theorem; the maps in question are simply linear, with no additional restrictions.  By contrast, the next result requires boundedness.  Let $\ooo$ denote the set of all outer functions in $\hdisk$.  
\begin{lem}[from {\cite[Lemma~4]{GiLa:Polya2011}}]\label{lem-Rank-1}
A bounded linear functional $\nu:\hdisk\rightarrow\complex$ satisfies
\[
\nu(\ooo)\subset\complex\setminus\{0\} \]
if and only if there exist a point $z_0\in\ddd$ and a scalar $\sigma\in\complex\setminus\{0\}$ such that for all $f\in \hdisk$,
\[
\nu(f)=\sigma f(z_0).
\]
\end{lem}

\section{The characterization problem\label{sec-Characterization}}

Theorem~\ref{thm-Structure} makes it possible to give describe explicitly the structure of a linear operator $A:\ltwo\rightarrow\ltwo$ for which $A(\tmm)\subset\tmm$.  
\begin{thm}\label{thm-Characterization}
Let $A:\ltwo\rightarrow\ltwo$ be a bounded linear operator such that $A(\tmm)\subset\tmm$.  Then there exist a function $\psi\in\sss(\ddd)\cap\hdisk$ and an analytic function $\varphi:\ddd\rightarrow\ddd$ such that 
\[
A=\hhh^{-1}M_\psi C_\varphi\hhh.
\]
\end{thm}
\begin{pf}
Consider the associated operator 
\[
B=\hhh A\hhh^{-1}:\hdisk\rightarrow\hdisk.
\]
Noting that every stable polynomial $p\in\poly(\ddd)$ is an outer function, the hypothesis that $A$ preserves translates of minimum phase functions implies that  
\[
A\hhh^{-1}p\in\tmm.
\]
Letting $g\in\tmm$ have the form $T_\tau f$, where $\tau\geq0$ and $f\in\mmm$, the definition of $\hhh$ yields that 
\[
(\hhh g)(z)=\exp\!\left(\tau\frac{z-1}{z+1}\right)(\hhh f)(z).
\]
Since $f$ is minimum phase, the function $\hhh f$ is outer.  But if $\tau>0$, then $\hhh g$ is not outer, since the function $e^{\tau\frac{z-1}{z+1}}$ is singular inner; see \cite[Section~2.6]{MaRo:2007}.   Nevertheless, $e^{\tau\frac{z-1}{z+1}}$ is $\ddd$-stable, as is every outer function in $\hdisk$, and therefore $\hhh g\in\sss(\ddd)$.  This shows that 
\[
B(\poly(\ddd))\subset\sss(\ddd),
\]
so that $B$ satisfies the hypothesis of Theorem~\ref{thm-Structure} with $\Omega=\ddd$.   Thus $B$ is either a product-composition operator or a rank one operator.  

By hypothesis $B$ does not map any outer function to 0.  This implies that if $B$ is a rank one operator of the form given in the first part of Theorem~\ref{thm-Structure}, then the associated linear functional $\nu$ satisfies the hypothesis of Lemma~\ref{lem-Rank-1}. The conclusion of Lemma~\ref{lem-Rank-1} then implies that $B$ has the form
\[
B=M_\psi C_\varphi,
\]
with $\varphi(z)=z_0$ a constant function.   This shows that $B$ may be represented as a product-composition operator whatever its rank, establishing the desired characterization of the original operator $A$.
\end{pf}

At this point the issue of sufficiency naturally arises, and there are two basic questions. What conditions on $\psi$ and $\varphi$ guarantee that:
\begin{enumerate}
\item $M_\psi C_\varphi:\hdisk\rightarrow\hdisk$ is a bounded operator;
\item the associated operator $A=\hhh^{-1}M_\psi C_\varphi\hhh$ maps $\tmm$ into itself?
\end{enumerate}
Regarding question~1, for any analytic function $\varphi:\ddd\rightarrow\ddd$ and any $\psi\in H^\infty_\ddd$, the product-composition operator $M_\psi C_\varphi$ is a bounded linear operator.  See \cite[Theorem~5.1.5]{MaRo:2007}.  But these conditions are not necessary in that $\psi$ need not be bounded.  See \cite[Section~4]{GiLaMa:JFA2011} for a more detailed discussion.   
Concerning question~2, the operator $A=\hhh^{-1}M_\psi C_\varphi\hhh$ preserves $\tmm$ if the product-composition operator $M_\psi C_\varphi$ preserves $\hhh(\tmm)$, which consists of functions of the form  
\[
\exp\!\left(\frac{\tau(z-1)}{z+1}\right)\rho(z),
\]
where $\rho\in\hdisk$ is outer and $\tau\geq0$.  If $\psi\in\hhh(\tmm)$ and either $\varphi(z)=z$ or the function
\[
\frac{1-\varphi(z)}{1+\varphi(z)} 
\]
has an analytic extension to the closed disk $\overline{\ddd}$, then $M_\psi C_\varphi$ is guaranteed to preserve $\hhh(\tmm)$.   These conditions are not particularly stringent, and the corresponding class of operators is large.  

Returning to the earlier thread, $A=\hhh^{-1}M_\psi C_\varphi\hhh$ has the integral representation 
\[
(Af)(t)=\frac{1}{\pi}\int_{y=-\infty}^\infty\int_{\tau=0}^\infty\eta(y)e^{\theta(y,\tau,t)}f(\tau)\,d\tau dy,
\]
where 
\[
\eta(y)=\frac{\frac{1}{1+iy}\psi\left(\frac{1-iy}{1+iy}\right)}{1+\varphi\left(\frac{1-iy}{1+iy}\right)}
\]
and 
\[
\theta(y,\tau,t)=\tau\frac{\varphi\left(\frac{1-iy}{1+iy}\right)-1}{\varphi\left(\frac{1-iy}{1+iy}\right)+1}+iyt.
\]
The above formulation is somewhat awkward on account of its representation in terms of functions on the unit disk.  An alternative is to express $A=\hhh^{-1}M_\psi C_\varphi\hhh$ in terms of functions on the right half plane, as follows.  
\begin{thm}\label{thm-Characterization2}
Let $A:\ltwo\rightarrow\ltwo$ be a continuous linear operator such that that $A(\tmm)\subset\tmm$.  Then there exist a function $\alpha\in\sss(\complex_+)\cap\hright$ and an analytic function $\xi:\complex_+\rightarrow\complex_+$ such that 
\[
A=\lll^{-1}M_\kappa C_\xi\lll,
\]
with $\kappa=\sqrt{2\pi}(1+\xi)\alpha$.  
\end{thm}
\begin{pf}
By Theorem~\ref{thm-Characterization} the operator $A$ has the form 
\[
A=\hhh^{-1}M_\psi C_\varphi\hhh=\lll^{-1}\Phi^{-1}M_\psi C_\varphi\Phi\lll,
\]
where $\psi\in\sss(\ddd)\cap\hdisk$ and $\varphi:\ddd\rightarrow\ddd$ is analytic.  The operator $\Phi^{-1}M_\psi C_\varphi\Phi$ is itself a product-composition operator, as follows.  For any $F\in\hright$,
\begin{eqnarray*}
(M_\psi C_\varphi\Phi F)(z)&=&\psi(z)\frac{2\sqrt{\pi}}{1+\varphi(z)}F\left(\frac{1-\varphi(z)}{1+\varphi(z)}\right)\\
&=&\sqrt{\pi}\psi(z)\frac{2}{1+\varphi(z)}(F\circ m\circ\varphi)(z)\\
&=&\sqrt{\pi}\psi(z)(1+m\circ\varphi)(F\circ m\circ\varphi)(z).
\end{eqnarray*}
Therefore
\begin{eqnarray*}
\Phi^{-1}M_\psi C_\varphi\Phi F&=&\sqrt{\pi}(\Phi^{-1}\psi)(1+m\circ\varphi\circ m)F\circ m\circ\varphi\circ m\\
&=& M_\kappa C_\xi F,
\end{eqnarray*}
where 
\[
\kappa=\sqrt{\pi}(\Phi^{-1}\psi)(1+m\circ\varphi\circ m)
\rule{5pt}{0pt}\mbox{ and }\rule{5pt}{0pt}\xi=m\circ\varphi\circ m.
\]  
Since $\varphi:\ddd\rightarrow\ddd$ is analytic, it follows that $\xi:\complex_+\rightarrow\complex_+$ is analytic too.  Note that $\Phi^{-1}$ carries $\ddd$-stable functions in $\hdisk$ to $\complex_+$-stable functions in $\hright$, and set 
\[
\alpha=\frac{1}{\sqrt{2}}\Phi^{-1}\psi,
\]
so that $\alpha\in\hright$ is $\complex_+$-stable, just as $\psi$ is $\ddd$-stable, the factor of $1/\sqrt{2}$ being for later convenience.  In terms of this notation, 
\[
A=\lll^{-1}M_\kappa C_\xi\lll\rule{5pt}{0pt}\mbox{ where }\rule{5pt}{0pt}\kappa=\sqrt{2\pi}(1+\xi)\alpha,
\]
completing the proof.  
\end{pf}

The formulation of Theorem~\ref{thm-Characterization2} is preferable to that of Theorem~\ref{thm-Characterization}, since an operator $A=\lll^{-1}M_\kappa C_\xi\lll$ with $\kappa=\sqrt{2\pi}(1+\xi)\alpha$ has an integral representation of the comparatively simple form
\begin{equation}\label{IntegralForm}
(Af)(t)=\frac{1}{\sqrt{2\pi}}\int_{y=-\infty}^\infty\int_{\tau=0}^\infty\alpha(iy)(1+\xi(iy))e^{-\tau\xi(iy)+iyt}f(\tau)\,d\tau dy.
\end{equation}
The reason for expressing $\kappa=\sqrt{2\pi}(1+\xi)\alpha$ in factored form is that, from the point of view of test data, $\alpha$ is the primary object.  In the next section it will be seen that $\lll^{-1}\alpha$ can be measured directly by evaluation of $A$ at the test function $e^{-t}$, whereas $\kappa$ must be computed from a combination of measurements.

\section{The identification problem\label{sec-Identification}}
With Theorem~\ref{thm-Characterization} in hand, the problem of identifying a given minimum phase preserving operator using test functions becomes straightforward.   
\begin{thm}\label{thm-Identification}
Let $A:\ltwo\rightarrow\ltwo$ have the property that $A(\tmm)\subset\tmm$.  Then $A$ can be reconstructed from its values at $\rho_0(t)=\sqrt{2}e^{-t}$ and $\rho_1(t)=\sqrt{2}e^{-t}(2t-1)$ by means of the formula $A=\hhh^{-1}M_\psi C_\varphi\hhh$, where 
\[
\psi=\hhh A\rho_0\rule{5pt}{0pt}\mbox{ and }\rule{5pt}{0pt}\varphi=(\hhh A\rho_1)/(\hhh A\rho_0).
\]
\end{thm}
\begin{pf}
Recall that 
\[
\rho_0(t)=(\hhh^{-1}\mu_0)(t)\rule{5pt}{0pt}\mbox{ and }\rule{5pt}{0pt}\rho_1(t)= (\hhh^{-1}\mu_1)(t),
\]
so that 
\[
A\rho_0=\hhh^{-1}\psi\rule{5pt}{0pt}\mbox{ and }\rule{5pt}{0pt}A\rho_1=\hhh^{-1}(\psi\varphi).
\]
\end{pf}
\begin{cor}\label{cor-Identification}
The function $\psi$ in the conclusion of Theorem~\ref{thm-Identification} has the form
\[
\psi(z)=\exp\!\left(\frac{\tau(z-1)}{z+1}\right)\psi_0(z),
\]
where $\tau\geq0$ and $\psi_0\in\hdisk$ is outer. 
\end{cor}
\begin{pf}
The test function $\rho_0$ is minimum phase, so by hypothesis $A\rho_0$ is the translate by some $\tau\geq0$ of a minimum phase function. It follows that $\hhh A\rho_0$ has the stated form, just as in the proof of Theorem~\ref{thm-Characterization}.  
\end{pf}

Note that the functions $\rho_0$ and $\rho_1$ in Theorem~\ref{thm-Identification} can be replaced by any pair having the same span.  The following choice is geared to fit Theorem~\ref{thm-Characterization2} and the formulation (\ref{IntegralForm}), resulting in simple reconstruction formulas.    
\begin{thm}\label{thm-Reformulation}
Every continuous linear operator $A:\ltwo\rightarrow\ltwo$ with the property that $A(\tmm)\subset\tmm$ is determined by its values $A\sigma_0$ and $A\sigma_1$, where
\[
\sigma_0(t)=e^{-t}(1-t)\rule{5pt}{0pt}\mbox{ and }\rule{5pt}{0pt}\sigma_1(t)=e^{-t}t
\]
for all $t\in\real_+$.  Setting $\alpha=\lll A\sigma_0+\lll A\sigma_1$ and $\xi={\lll A\sigma_0}/{\lll A\sigma_1}$,
\[
A=\lll^{-1}M_\kappa C_\xi\lll \rule{5pt}{0pt}\mbox{ with }\rule{5pt}{0pt}\kappa=\sqrt{2\pi}(1+\xi)\alpha;
\]
equivalently, $A$ has the integral representaion
\[
(Af)(t)=\frac{1}{\sqrt{2\pi}}\int_{y=-\infty}^\infty\int_{\tau=0}^\infty\alpha(iy)(1+\xi(iy))e^{-\tau\xi(iy)+iyt}f(\tau)\,d\tau dy.
\]
Moreover, $\alpha$ necessarily has the form 
$
\alpha(w)=e^{-\varepsilon w}\alpha_0(w)
$
for some scalar $\varepsilon\geq0$ and some outer function $\alpha_0\in\hright$, and $\xi:\complex_+\rightarrow\complex_+$ is analytic.  
\end{thm}
\begin{pf}
Note that
\begin{eqnarray*}
(\lll \sigma_0)(w)&=&\frac{w}{\sqrt{2\pi}(1+w)^2},\\
(\lll \sigma_1)(w)&=&\frac{1}{\sqrt{2\pi}(1+w)^2}.
\end{eqnarray*}
Using these formulas and Theorem~\ref{thm-Characterization2}, a straightforward calculation gives the desired formulas for $\alpha$ and $\xi$.   Since $\sigma_0+\sigma_1$ is minimum phase, the hypothesis that $A$ preserves translates of minimum phase functions implies that $\alpha$ is the Fourier-Laplace transform of a translate by some $\varepsilon\geq0$ of a minimum phase function.  It follows that $\alpha(w)=e^{-\varepsilon w}\alpha_0(w)$ for some outer function $\alpha_0\in\hright$.  
\end{pf}

There is no way to identify an unknown operator of the form $A=\hhh^{-1}M_\psi C_\varphi\hhh$ from a single value $Af$, no matter which test function $f\in\ltwo$ is chosen.  In this sense Theorems~\ref{thm-Identification} and \ref{thm-Reformulation} give results that are best possible.  

%(Af)(t)=\frac{1}{\sqrt{2\pi}}\int_{y=-\infty}^\infty\int_{\tau=0}^\infty\alpha(iy)(1+\xi(iy))e^{-\tau\xi(iy)+iyt}f(\tau)\,d\tau dy

\section{Remarks\label{sec-Remarks}}

\subsection{Operators preserving $\mmm$\label{sec-Operators}}
The identification of bounded linear operators $A:\ltwo\rightarrow\ltwo$ preserving minimum phase functions, but not necessarily their translates, can be handled using exactly the same methods as above.  Indeed without considering translates, the arguments become simpler, and the analogue of Theorem~\ref{thm-Reformulation} that results is only slightly different, as follows. 
\begin{thm}\label{thm-Minimum}
Every continuous linear operator $A:\ltwo\rightarrow\ltwo$ with the property that $A(\mmm)\subset\mmm$ is determined by its values $A\sigma_0$ and $A\sigma_1$, where
\[
\sigma_0(t)=e^{-t}(1-t)\rule{5pt}{0pt}\mbox{ and }\rule{5pt}{0pt}\sigma_1(t)=e^{-t}t
\]
for all $t\in\real_+$.  Setting $\alpha=\lll A\sigma_0+\lll A\sigma_1$ and $\xi={\lll A\sigma_0}/{\lll A\sigma_1}$,
\[
A=\lll^{-1}M_\kappa C_\xi\lll \rule{5pt}{0pt}\mbox{ with }\rule{5pt}{0pt}\kappa=\sqrt{2\pi}(1+\xi)\alpha;
\]
equivalently, $A$ has the integral representaion
\[
(Af)(t)=\frac{1}{\sqrt{2\pi}}\int_{y=-\infty}^\infty\int_{\tau=0}^\infty\alpha(iy)(1+\xi(iy))e^{-\tau\xi(iy)+iyt}f(\tau)\,d\tau dy.
\]
Moreover, $\alpha\in\hright$ is necessarily an outer function, and $\xi:\complex_+\rightarrow\complex_+$ is analytic.  
\end{thm}
Thus the only change is that the function $\alpha$ in the integration kernel is necessarily outer, whereas before it was allowed to be a modulated outer function.   Another way to approach the case $A(\mmm)\subset\mmm$ is to use the minimum phase preserving isometry
\[
\mathfrak{d}:\ltwo\rightarrow\lltwo
\]
defined in Section~\ref{sec-Isometry}, and appeal to the corresponding theorem for the discrete case, \cite[Theorem~4]{GiLa:Polya2011}.  The mapping $\mathfrak{d}$ does not carry translates of continuous signals to translates of discrete signals, so by contrast Theorem~\ref{thm-Reformulation} cannot be proved using its discrete analogue.  Indeed, whereas Theorem~\ref{thm-Reformulation} needs the full force of Theorem~\ref{thm-Structure}, its discrete analogue was first proved in \cite[Theorem~4]{GiLaMa:JFA2011} using less powerful methods.

\subsection{The hypothesis of boundedness\label{sec-Hypothesis}}

It is a remarkable fact that if the operator $A:\ltwo\rightarrow\ltwo$ is assumed to have rank at least two, then the hypothesis of boundedness can be dropped from Theorems~\ref{thm-Characterization}, \ref{thm-Characterization2}, \ref{thm-Identification}, \ref{thm-Reformulation} and \ref{thm-Minimum}.  Boundedness was only used when invoking  Lemma~\ref{lem-Rank-1}, which applies to the rank 1 case.  Moreover, in the case of rank at least two, the property that $A(\tmm)\subset\tmm$ or $A(\mmm)\subset\mmm$ forces $A$ to be injective, a property implied by the product-composition structure of $\hhh A\hhh^{-1}$.  See \cite[Corollary~2]{GiLa:Polya2011}.

\section{Conclusion\label{sec-Conclusion}}

This paper's main results, Theorems~\ref{thm-Reformulation} and \ref{thm-Minimum} in Sections~\ref{sec-Identification} and \ref{sec-Operators}, solve completely the natural inverse problem concerning linear operators on the half line that preserve the class of minimum phase---or translated minimum phase---signals.  They show that every such operator is conjugate via the Fourier-Laplace transform to a product-composition operator on Hardy space, that precisely two test functions are required to reconstruct such an operator, and that there is a simple reconstruction formula in terms of the exponentially damped linear test functions $e^{-t}(1-t)$ and $e^{-t}t$.    Recent developments in the theory of stable polynomials played a crucial role in the derivation of these results.   

The inverse problem of identifying a minimum phase preserving operator from its output was motivated by considerations in seismic imaging, where seismic traces are widely assumed to be minimum phase functions, or translates of minimum phase functions.   Results in the present paper show that if the minimum phase hypothesis is valid for seismic data, then current techniques such as Wiener deconvolution are inadequate in the most general setting.  A non-stationary operator representing the process of an input wavelet travelling through an interval of rock, if it preserves minimum phase, cannot be determined from a single test function.  And therefore an input wavelet cannot be computed from a single seismic trace---even with non-stationary methods; at least two traces, arising from two different source functions are required to determine the operator, and thus to invert it.   This offers the exciting possibility that a re-evaluation of the seismic deconvolution techniques in light of the Theorem~\ref{thm-Reformulation} could lead to substantial improvement.  In any case the very precise results in the present paper shine a spotlight on the minimum phase hypothesis itself, and call for a better understanding of its origins in models of layered media.

\bibliography{ReferencesOuter}
\end{document}